\documentclass[a4paper,12pt]{amsart}
\usepackage{amsfonts}

\usepackage{amsmath,latexsym,amscd}
\usepackage{amssymb} 
\usepackage[mathscr]{eucal} 
\usepackage[all]{xy}

\title{Complex projective threefolds with non-negative canonical
Euler-Poincare characteristic}
\author{Meng Chen and Kang Zuo}
\address{\rm Institute of Mathematics, School of Mathematical
Sciences, Fudan University, Shanghai, 200433, China}
\email{mchen@fudan.edu.cn}
\address{\rm Faculty of Mathematics \& Computer Science,
Johannes Gutenberg Universitaet Mainz, Staudingerweg 9, D-55128
Mainz, Germany}
\email{kzuo@mathematik.uni-mainz.de}
\thanks{This paper was
supported by both the Program for New Century Excellent
Talents in University (\#NCET-05-0358) and the National
Outstanding Young Scientist Foundation (\#10625103)}
\newcommand{\bC}{{\mathbb C}}
\newcommand{\bQ}{{\mathbb Q}}

\newcommand{\roundup}[1]{\ulcorner{#1}\urcorner}

\newtheorem{thm}{Theorem}[section]
\newtheorem{lem}[thm]{Lemma}
\newtheorem{cor}[thm]{Corollary}
\newtheorem{prop}[thm]{Proposition}
\newtheorem{claim}[thm]{Claim}
\theoremstyle{definition}
\newtheorem{defn}[thm]{Definition}
\newtheorem{setup}[thm]{}

\newtheorem{exmp}[thm]{Example}

\newtheorem{rem}[thm]{Remark}
\theoremstyle{remark}

\def\eea{\end{eqnarray*}}
\def\bea{\begin{eqnarray*}}

\begin{document}
\begin{abstract}
Let $V$ be a complex nonsingular projective 3-fold of
general type with $\chi(\omega_V)\geq 0$ (resp. $>0$). We
prove that the m-canonical map $\Phi_{|mK_V|}$ is
birational onto its image for all $m\ge 14$ (resp. $\geq
8$). Known examples show that the lower bound $r_3=14$
(resp. $=8$) is optimal.
\end{abstract}
\maketitle
\pagestyle{myheadings} \markboth{\hfill M. Chen and K. Zuo
\hfill}{\hfill  \hfill}
\section{\bf Introduction}
We work over the complex number field $\bC$.

In this paper we study multi-canonical systems $|mK_V|$ on complex
projective threefolds $V$ of general type.

On a smooth complex complete curve $C$ of genus $g(C)\ge 2$ it is
well-known that the 3-canonical map $\varphi_3:=\Phi_{|3K_C|}$ is
always an embedding. On a smooth minimal projective surface $S$ of
general type Bombieri \cite{Bom} proved that the m-canonical map is
a birational morphism for all $m\ge 5$. The work of Tsuji
\cite{Tsuji}, Hacon-M$^{\rm c}$Kernan \cite{H-M} and Takayama
\cite{Tak} says that there exists a universal constant $r_3$ such
that the $r_3$-canonical map $\Phi_{r_3}$ is birational for all
smooth projective 3-folds of general type. We note that Tsuji
\cite{Tsuji} announced a very large $r_3$. Recently J.A. Chen and
the first author \cite{Jungkai-Meng} have given an explicit $r_3$
($\leq 77$) .

Under extra assumptions there have been already the
following optimal results about $r_3$ for minimal
projective 3-folds $X$ of general type with
$\mathbb{Q}$-factorial terminal singularities:
\begin{quote}
\item $\bullet$ $\varphi_5:=\Phi_{|5K_X|}$ is birational if either $X$ is
Gorenstein (by J.A. Chen, M.Chen and D-Q.Zhang \cite{Crelle}) or
$p_g(X)\ge 4$ (by M. Chen \cite{IJM}) or $K_X^3\gg 0$ (by G. T.
Todorov \cite{Tod}). $r_3=5$ is optimal.

\item $\bullet$ $\varphi_8:=\Phi_{|8K_X|}$ is birational if $p_g(X)\ge 2$
(by M. Chen \cite{IJM}). $r_3=8$ is optimal.

\item $\bullet$ $\varphi_7$ is birational if $q(X)>0$
(by J. A. Chen and C. D. Hacon \cite{C-H}).

\item $\bullet$ $\varphi_5$
is birational if $\chi(\omega_X)\geq 0$ and $q(X)>0$ (by J.
A. Chen and C. D. Hacon (\cite{C-H06}). $r_3=5$ is optimal.
\end{quote}
\bigskip

 Going on the study of 3-folds with $\chi(\omega)\geq 0$, we
prove the following:

\begin{thm}\label{main} Let $V$ be a nonsingular projective 3-fold of general
type with
$\chi(\omega_V):=\sum_{i=0}^3(-1)^ih^i(V,\omega_V)\ge 0$
where $\omega_V$ is the canonical line bundle of $V$. Then
the m-canonical map $\varphi_m$ is birational onto its
image for all $m\geq 14$.
\end{thm}

The following example of A. R. Iano-Fletcher shows that $r_3=14$ in
Theorem \ref{main} is optimal.

\begin{exmp}\label{examp1} (See No.19 at page 151 in \cite{C-R}) The canonical
hypersurface $X_{28}\subset \mathbb{P}(1,3,4,5,14)$ has 3
terminal quotient singularities, $p_g(X)=1$,
$q(X)=h^2(\mathcal{O}_X)=0$ and
$\chi(\omega_X)=-\chi(\mathcal {O}_X)=0$. A smooth model
$V$ of $X_{28}$ has the invariant: $\chi(\omega_V)=0$ and
$V$ is of general type. It is clear that $\varphi_m$ is
birational for all $m\ge 14$ and that $\varphi_{13}$ is not
birational. So the lower bound $r_3=14$ in Theorem
\ref{main} is sharp.
\end{exmp}

Our method has a direct consequence:

\begin{cor}\label{corollary} Let $V$ be a smooth projective 3-fold of
general type with $\chi(\omega_V)> 0$. Then the m-canonical
map $\Phi_{|mK_V|}$ is birational onto its image for all
$m\ge 8$.
\end{cor}

\begin{exmp} (See No.12 at page 151 in \cite{C-R}) Fletcher has a 3-fold
$X_{16}\subset \mathbb{P}(1,1,2,3,8)$ on which the
8-canonical map is birational and the 7-canonical map is
not birational. A smooth model $V$ of $X_{16}$ has the
invariant $\chi(\omega)=1$ and $V$ is of general type. $V$
has $r_3(V)=8$. Thus the lower bound $r_3=8$ in Corollary
\ref{corollary} is again optimal.
\end{exmp}

Note that 3-folds of general type with $\chi(\omega)\ge 0$
form an interesting class, around which there have been
already some established works:
\begin{quote}
\item$\bullet$ Gorenstein minimal 3-folds of general type have
$\chi(\omega)>0$;

\item$\bullet$ Any 3-fold of general type admitting a generically finite
cover over an Abelian variety has $\chi(\omega)>0$ (see
Green-Lazarsfeld \cite{G-L});

\item$\bullet$ Any 3-fold of general type, admitting a representation
$\rho: \pi_1(X)\rightarrow \text{GL}$ such that the Shafarevich map
$\text{Sh}_\rho$ is generically finite, has $\chi(\omega)\ge 0$ (see
Jost-Zuo \cite{J-Z});

\item$\bullet$ Any 3-fold of general type with generically large
fundamental group is conjectured to have $\chi(\omega)>0$ (see
Koll\'ar \cite{Kollar}, Conjecture 18.12.1).

\end{quote}
\medskip

It is clear that the key point in studying pluricanonical maps is
to compute $P_m$ for $m>1$. One may do this on a minimal model
according to Reid \cite{YPG}. When $\chi(\omega)=-\chi({\mathcal
O})\geq 0$, it is clear from Reid's formula that $P_2>0$. This
makes it possible for us to prove effective birationality of
$\varphi_m$. But to prove the optimal bound $r_3=14$ need more
evolved argument, which is exactly what we have done in this paper.
The case $\chi(\omega)<0$ has been treated with quite a different
approach very recently by J. A. Chen and the first author
\cite{Jungkai-Meng}.

\section{\bf Notations and set up}

Since both the birationality of pluricanonical maps and
$\chi(\omega)$ are birationally invariant we may study a
minimal model. By the 3-dimensional MMP (see \cite{K-M,
KMM} for instance) we only have to consider a minimal
3-fold $X$ of general type with $\mathbb{Q}$-factorial
terminal singularities. Denote the Cartier index of $X$ by
$r:=r(X)$ which is the minimal positive integer with $rK_X$
a Cartier divisor, where $K_X$ is a canonical divisor on
$X$. The symbol $\equiv$ stands for the numerical
equivalence of divisors, whereas $\sim$ denotes the linear
equivalence and $=_{\mathbb{Q}}$ denotes the
$\mathbb{Q}$-linear equivalence.

\begin{setup}\label{setup}{\bf Set up for $|m_0K_X|$.}
Assume $P_{m_0}(X):=h^0(X, \mathcal {O}_X(m_0K_X))\ge 2$ for some
positive integer $m_0>0.$ We study the $m_0$-canonical map
$\varphi_{m_0}$ which is a rational map.

First we fix an effective Weil divisor $K_{m_0}\sim m_0K_X$. Take
successive blow-ups $\pi: X'\rightarrow X$ (along nonsingular
centers), which exists by Hironaka's big theorem, such that:
\par

(i) $X'$ is smooth;

(ii) the movable part of $|m_0K_{X'}|$ is base point free;

(iii) the support of $\pi^*(K_{m_0})$ is of simple normal crossings.
\par

Denote by $g_{m_0}$ the composition $\varphi_{m_0}\circ\pi$. So
$g_{m_0}: X'\longrightarrow W'\subseteq{\mathbb P}^{P_{m_0}(X)-1}$
is a morphism. Let $X'\overset{f_{m_0}}\longrightarrow
B\overset{s}\longrightarrow W'$ be the Stein factorization of
$g_{m_0}$. We have the following commutative diagram:\medskip

\begin{picture}(50,80) \put(100,0){$X$} \put(100,60){$X'$}
\put(170,0){$W'$} \put(170,60){$B$} \put(112,65){\vector(1,0){53}}
\put(106,55){\vector(0,-1){41}} \put(175,55){\vector(0,-1){43}}
\put(114,58){\vector(1,-1){49}} \multiput(112,2.6)(5,0){11}{-}
\put(162,5){\vector(1,0){4}} \put(133,70){$f_{m_0}$}
\put(180,30){$s$} \put(92,30){$\pi$}
\put(135,-8){$\varphi_{m_0}$}\put(136,40){$g_{m_0}$}
\end{picture}
\bigskip

We recall the definition of $\pi^*(K_X)$. One has
$r(X)K_{X'}=\pi^*(r(X)K_X)+E_{\pi}$ where $E_{\pi}$ is a sum of
effective exceptional divisors. One defines
$\pi^*(K_X):=K_{X'}-\frac{1}{r(X)}E_{\pi}$. So, whenever we take
the round up of $m\pi^*(K_X)$, we always have
$\roundup{m\pi^*(K_X)}\leq mK_{X'}$ for any integer $m>0$. We may
write $m_0\pi^*(K_X)=_{\mathbb Q} M_{m_0}+E_{m_0}',$ where
$E_{m_0}'$ is an effective  $\bQ$-divisor and $M_{m_0}$ is the
movable part of $|m_0K_{X'}|$. On the other hand, one has
$m_0K_{X'}=_{\mathbb Q}\pi^*(m_0K_X)+E_{\pi, m_0}=M_{m_0}+Z_{m_0},$
where $Z_{m_0}$ is the fixed part and $E_{\pi, m_0}$ an effective
${\mathbb Q}$-divisor which is a ${\mathbb Q}$-sum of distinct
exceptional divisors. Clearly $Z_{m_0}=E_{m_0}'+E_{\pi, m_0}$.

If $\dim\varphi_{m_0}(X)=2$, a general fiber of $f_{m_0}$ is a
smooth projective curve of genus $\ge 2$. We say that $X$ is {\it
$m_0$-canonically fibred by curves}.

If $\dim\varphi_{m_0}(X)=1$, a general fiber $S$ of
$f_{m_0}$ is a smooth projective surface of general type.
We say that $X$ is {\it $m_0$-canonically fibred by
surfaces} with invariants $(c_1^2(S_0), p_g(S)),$ where
$S_0$ is the minimal model of $S$. We may write
$M_{m_0}\equiv a_{m_0}S$ where $a_{m_0}\ge P_{m_0}(X)-1$ by
considering the degree of a curve in a projective space.

{\it A generic irreducible element $S$ of} $|M_{m_0}|$ means either
a general member of $|M_{m_0}|$ whenever $\dim\varphi_{m_0}(X)\ge 2$
or, otherwise, a general fiber of $f_{m_0}$.
\end{setup}

\begin{defn}\label{generic} Assume that $|M'|$ is movable on $V$.
By abuse of concepts, we also define {\it a generic
irreducible element} $S'$ of an arbitrary linear system
$|M'|$ on an arbitrary variety $V$ in a similar way. {\it A
generic irreducible element} $S'$ of $|M'|$ is defined to
be a generic irreducible component in a general member of
$|M'|$.
\end{defn}

\section{\bf A technical theorem}

Believing that Theorem 2.2 in \cite{IJM} is quite effective
in treating 3-folds $X$ with $p_g(X)\ge 2$, we extend the
technique there to build a parallel theorem so as to study
those $X$ with $p_g(X)\le 1$, but with $P_{m_0}(X)\ge 2$
for some integer $m_0>0$.

\begin{setup}\label{a}{\bf Assumptions}. We need to make
the following assumptions to explain our key method. Keep
the same notation as in $\ref{setup}$ above. Let $m>0$ be
certain integer:
\begin{itemize}

\item[(i)] Either $m\geq
m_0+2$ and $p_g(X)>0$ or $|mK_{X'}|$ separates different
irreducible elements $S$ of $|M_{m_0}|$ (namely,
$\Phi_{|mK_{X'}|}(S')\neq \Phi_{|mK_{X'}|}(S'')$ for two
different irreducible elements $S'$, $S''$ of $|M_{m_0}|$)
and $p_g(X)=0$.

\item[(ii)] Assume
that, on the smooth surface $S$, there is a movable linear
system $|G|$ and that $C$,as a generic irreducible element
of $|G|$, is smooth. The linear system $|mK_{X'}||_S$ on
$S$ (as a sub-linear system of $|mK_{X'}|_S|$) separates
different generic irreducible elements of $|G|$. Or
sufficiently, the complete linear system
$$|K_{S}
+\roundup{(m-1)\pi^*(K_X)-S-\frac{1}{p}E_{m_0}'}|_{S}|$$
separates different generic irreducible elements of $|G|$.

\item[(iii)]
There is a rational number $\beta>0$ such that
${\pi}^*(K_X)|_{S}-\beta C$ is numerically equivalent to an
effective ${\mathbb Q}$-divisor; {\bf Set}
$$\alpha:=(m-1-\frac{m_0}{p}-\frac{1}{\beta})\xi$$ and
$\alpha_0:=\roundup{\alpha}$.
\item[(iv)]
Either the inequality $\alpha > 1$ holds; or $C$ is
non-hyperelliptic, $m-1-\frac{m_0}{p}-\frac{1}{\beta}>0$
and $C$ is an even divisor on $S$.
\item[(v)]
Either $\alpha>2$; or $\alpha_0\ge 2$  and $C$ is
non-hyperelliptic; or $C$ is non-hyperelliptic,
$m-1-\frac{m_0}{p}-\frac{1}{\beta}>0$ and $C$ is an even
divisor on $S$.
\end{itemize}
\end{setup}

Set $\xi:=(\pi^*(K_X)\cdot C)_{X'}$ which is a positive
rational number and define
$$p:=\begin{cases} 1  &\text{if}\ \dim (B)\ge 2\\
a_{m_0} (\text{see}\ \ \ \ref{setup}\ \text{for the
definition}) &\text{otherwise.}
\end{cases}$$
Let $f:=f_{m_0}: X'\longrightarrow B$ be an induced
fibration by $\varphi_{m_0}$.

\begin{thm}\label{Key} Let $X$ be a minimal
projective $3$-fold of general type with $P_{m_0}(X)\ge 2$
for some integer $m_0>0$. Keep the same notation as in
$\ref{setup}$ above. Then the inequality
$$m\xi\ge
2g(C)-2+\alpha_0$$ holds under Assumption \ref{a} $(iii)$
and $(iv)$. Furthermore $\varphi_m$ of $X$ is birational
onto its image under Assumption \ref{a} $(i), (ii), (iii)$
and $(v)$.
\end{thm}

\begin{proof} First we shall show that $|mK_{X'}|$
can separate different irreducible elements of $|M_{m_0}|$. When
$p_g(X)=0$, $|mK_{X'}|$ can separate different irreducible elements
of $|M_{m_0}|$ by assumption (i). When $p_g(X)>0$, we consider the
sub-system
$$|K_{X'}+\roundup{(m-1)\pi^*(K_X)-\frac{1}{p}E_{m_0}'}|\subset |mK_{X'}|.$$
Let $S'$ and $S''$ be two different generic irreducible elements of
$|M_{m_0}|$. Clearly one has
\begin{eqnarray*}
&&K_{X'}+\roundup{(m-1)\pi^*(K_X)-\frac{1}{p}E_{m_0}'}\\
&\ge& K_{X'}+\roundup{(m-m_0-1)\pi^*(K_X)}+S\ge S
\end{eqnarray*}
by assumption (i). So
$|K_{X'}+\roundup{(m-1)\pi^*(K_X)-\frac{1}{p}E_{m_0}'}|$ can
separate $S'$ and $S''$ if either $\dim (B)\ge 2$ (cf. Lemma 2 of
\cite{T}) or $\dim (B)=1$ and $g(B)=0$ (cf. 2.1(P2) of
\cite{MPCPS}). For the case $\dim (B)=1$ and $g(B)>0$, one has
$a_{m_0}\ge P_{m_0}\ge 2$. Thus $p\ge 2$. Since
$$m-1-\frac{2m_0}{p}\ge 1+(1-\frac{2}{p})m_0>0$$
and then
$$(m-1)\pi^*(K_X)-\frac{2}{p}E_{m_0}'-S'-S''\equiv
(m-1-\frac{2m_0}{p})\pi^*(K_X)$$ is nef and big, the
Kawamata-Viehweg vanishing theorem (\cite{Kav, V}) gives a
surjective map:
\begin{eqnarray*}
&& H^0(X',K_{X'}+\roundup{(m-1)\pi^*(K_X)-\frac{2}{p}E_{m_0}'})\\
&\longrightarrow& H^0(S',
K_{S'}+\roundup{(m-1)\pi^*(K_X)-\frac{2}{p}E_{m_0}'}|_{S'})\oplus\\
&&H^0(S'',
K_{S''}+\roundup{(m-1)\pi^*(K_X)-\frac{2}{p}E_{m_0}'}|_{S''}).
\end{eqnarray*}
The last two groups are non-zero because $p_g(X)>0$ (so $p_g(S')$,
$p_g(S'')>0$) and $m-1-m_0> 0$. Therefore
$|K_{X'}+\roundup{(m-1)\pi^*(K_X)-\frac{2}{p}E_{m_0}'}|$ can
separate $S'$ and $S''$ and so can $|mK_{X'}|$.

 By the birationality principle (P1)
and (P2) of \cite{MPCPS}, it suffices to prove that $|mK_{X'}||_S$
on $S$ gives a birational map onto its image. Practically we may
study a smaller linear system than $|mK_{X'}||_S$ on $S$. Noting
that $(m-1)\pi^*(K_X)-\frac{1}{p}E_{m_0}'-S\equiv
(m-1-\frac{m_0}{p})\pi^*(K_X)$ is nef and big under the assumptions
(iv) or (v), the vanishing theorem gives a surjective map
$$H^0(X',K_{X'}+\roundup{(m-1)\pi^*(K_X)-\frac{1}{p}E_{m_0}'})$$
$$\longrightarrow  H^0(S,
K_{S}+\roundup{(m-1)\pi^*(K_X)-S-\frac{1}{p}E_{m_0}'}|_{S}).\eqno
(3.1)$$ Note that
$|K_{X'}+\roundup{(m-1)\pi^*(K_X)-\frac{1}{p}E_{m_0}'}|\subset
|mK_{X'}|$. It suffices to prove that
$|K_{S}+\roundup{(m-1)\pi^*(K_X)-S-\frac{1}{p}E_{m_0}'}|_{S}|$ gives
a birational map.

The birationality principle again allows us to study the restriction
to curves by assumption (ii). Now consider a generic irreducible
element $C\in |G|$. By assumption (iii), there is an effective
${\mathbb Q}$-divisor $H$ on $S$ such that
$$\frac{1}{\beta}\pi^*(K_X)|_{S}\equiv C+H.$$
By the vanishing theorem, whenever
$m-1-\frac{m_0}{p}-\frac{1}{\beta}>0$, we have the surjective map
$$H^0(S,
K_{S}+\roundup{((m-1)\pi^*(K_X)-S-\frac{1}{p}E_{m_0}')|_{S}-H})$$
$$\longrightarrow  H^0(C, K_C + D)\hskip6.8cm \eqno (3.2)$$
where
$D:=\roundup{((m-1)\pi^*(K_X)-S-\frac{1}{p}E_{m_0}')|_{S}-C-H}|_C$
is a divisor on $C$. Noting that
$$((m-1)\pi^*(K_X)-S-\frac{1}{p}E_{m_0}')|_{S}-C-H\equiv (m-1-
\frac{m_0}{p}-\frac{1}{\beta})\pi^*(K_X)|_{S}$$ and that $C$ is nef
on $S$, we have $\deg(D)\geq\alpha$ and thus $\deg(D)\geq \alpha_0$.
Whenever $C$ is non-hyperelliptic,
$m-1-\frac{m_0}{p}-\frac{1}{\beta}>0$ and $C$ is an even divisor on
$S$, $\deg(D)\geq 2$ automatically follows and thus $|K_C+D|$ gives
a birational map. Whenever $\deg(D)\ge 3$, then
$$|K_{S}+\roundup{((m-1)\pi^*(K_X)-S-\frac{1}{p}E_{m_0}')|_{S}-H}||_C$$
gives a birational map. Since
\begin{eqnarray*}
&&|K_{S}+\roundup{((m-1)\pi^*(K_X)-S-\frac{1}{p}E_{m_0}')|_{S}-H}|\\
&\subset
&|K_{S}+\roundup{(m-1)\pi^*(K_X)-S-\frac{1}{p}E_{m_0}'}|_{S}|,
\end{eqnarray*}
the latter linear system gives a birational map. So $\varphi_m$ of
$X$ is birational.

Finally we show the inequality for $\xi$. Whenever we have
$\deg(D)\ge 2$, $|K_C+D|$ is base point free by the curve theory.
Denote by $|M_m|$ the movable part of $|mK_{X'}|$ and by $|N_m|$ the
movable part of
$|K_{S}+\roundup{((m-1)\pi^*(K_X)-S-\frac{1}{p}E_{m_0}')|_{S}-H}|$.
Applying Lemma 2.7 of \cite{MPCPS} to surjective maps (3.1) and
(3.2), we have
$$m\pi^*(K_X)|_{S}\ge N_m\ \ \text{and}\ \ (N_m\cdot C)_{S}\ge 2g(C)-2
+\deg (D).$$ Note that the above inequality holds without conditions
(i) and (ii). We are done.
\end{proof}

\begin{rem}\label{remark} If we replace $M_{m_0}$ in Theorem \ref{Key}
by any divisor $N_{m_0}\le M_{m_0}$ with $h^0(X', N_{m_0})\ge 2$,
Theorem \ref{Key} is still true accordingly. This is clear by the
proof. The main point is that it suffices to prove that a sub-linear
system of $|mK_{X'}|$ gives a birational map. To avoid frustrating
setting up and more complicated notations, we omit the proof in
details. The idea is, however, trivially similar.
\end{rem}

While applying Theorem \ref{Key}, one has to choose a suitable
movable system $|G|$ on $S$. Then quite a technical problem is to
find a suitable $\beta$ as in Theorem \ref{Key}(iii). The following
lemma presents the way for the most difficult case - the rational
pencil case.

\begin{lem}\label{beta} Keep the same notation as in \ref{setup} and
Theorem \ref{Key}. Assume $B=\mathbb{P}^1$. Let $f:X'\longrightarrow
\mathbb{P}^1$ be an induced fibration of $\varphi_{m_0}$. Denote by
$F:=S$ a general fiber of $f$. Then one can find a sequence of
rational numbers $\{\beta_n\}$ with $\lim_{n\mapsto +\infty} \beta_n
= \frac{p}{m_0+p}$ such that
$\pi^*(K_X)|_F-\beta_n\sigma^*(K_{F_0})$ is ${\bQ}$-linearly
equivalent to an effective ${\bQ}$-divisor $N_n$, where
$\sigma:F\longrightarrow F_0$ is the blow down onto the smooth
minimal model.
\end{lem}

\begin{proof}

One has $\mathcal {O}_{B}(p)\hookrightarrow {f}_*\omega_{X'}^{m_0}$
and therefore ${f}_*\omega_{X'/B}^{t_0p}\hookrightarrow
{f}_*\omega_{X'}^{t_0p+2t_0m_0}$ for any big integer $t_0$.

For any positive integer $k$, denote by $M_k$ the movable part of
$|kK_{X'}|$. Note that ${f}_*\omega_{X'/B}^{t_0p}$ is generated by
global sections since it is semi-positive according to E. Viehweg
(\cite{VV}). So any local section can be extended to a global one.
On the other hand, $|t_0p\sigma^*(K_{F_0})|$ is base point free and
is exactly the movable part of $|t_0pK_F|$ by Bombieri \cite{Bom} or
Reider \cite{Reider}. Clearly one has the following relation:
$$a_0\pi^*(K_X)|_F\ge M_{t_0p+2t_0m_0}|_F\ge b_0\sigma^*(K_{F_0})$$
where $a_0:=t_0p+2t_0m_0$ and $b_0:=t_0p$. This means that there is
an effective $\mathbb{Q}$-divisor $E_0'$ on $F$ such that
$$a_0\pi^*(K_X)|_F=_{\bQ} b_0\sigma^*(K_{F_0})+E_0'.$$
Thus $\pi^*(K_X)|_F =_{\bQ} \frac{p}{p+2m_0}\sigma^*(K_{F_0})+E_0$
with $E_0=\frac{1}{a_0}E_0'$.

We consider the case $p\ge 2$.

Assume that we have defined $a_n$ and $b_n$ such that the following
is satisfied with $l = n:$
$$a_{l}\pi^*(K_X)|_F \ge b_{l}\sigma^*(K_{F_0}).$$
We will define $a_{n+1}$ and $b_{n+1}$ inductively such that the
above inequality is satisfied with $l = n+1$. One may assume from
the beginning that $a_n\pi^*(K_X)$ supports on a divisor with normal
crossings. Then the Kawamata-Viehweg vanishing theorem implies the
surjective map
$$H^0(K_{X'}+\roundup{a_n\pi^*(K_X)}+F)\longrightarrow H^0(F, K_F+
\roundup{a_n\pi^*(K_X)}|_F).$$
One has the relation
\begin{eqnarray*}
|K_{X'}+\roundup{a_n\pi^*(K_X)}+F||_F&=&|K_F+\roundup{a_n\pi^*(K_X)}|_F|\\
&\supset& |K_F+b_n\sigma^*(K_{F_0})|\\
&\supset& |(b_n+1)\sigma^*(K_{F_0})|.
\end{eqnarray*}
Denote by $M_{a_n+1}'$ the movable part of $|(a_n+1)K_{X'}+F|$.
Applying Lemma 2.7 of \cite{MPCPS} again, one has $M_{a_n+1}'|_F\ge
(b_n+1)\sigma^*(K_{F_0}).$ Re-modifying our original $\pi$ such that
$|M_{a_n+1}'|$ is base point free. In particular, $M_{a_n+1}'$ is
nef. Since $X$ is of general type $|mK_X|$ gives a birational map
whenever $m$ is big enough. Thus we see that $M_{a_n+1}'$ is big if
we fix a very big $t_0$ in advance.

Now the Kawamata-Viehweg vanishing theorem again gives
\begin{eqnarray*}
|K_{X'}+M_{a_n+1}'+F||_F&=&|K_F+M_{a_n+1}'|_F|\\
&\supset& |K_F+(b_n+1)\sigma^*(K_{F_0})|\\
&\supset& |(b_n+2)\sigma^*(K_{F_0})|.
\end{eqnarray*}

We may repeat the above procedure inductively. Denote by
$M_{a_n+t}'$ the movable part of $|K_{X'}+M_{a_n+t-1}'+F|$ for $t\ge
2$. For the same reason, we may assume $|M_{a_n+t}'|$ to be base
point free. Inductively one has:
$$M_{a_n+t}'|_F\ge (b_n+t)\sigma^*(K_{F_0}).$$
Applying the vanishing theorem once more, we have
\begin{eqnarray*}
|K_{X'}+M_{a_n+t}'+F||_F&=&|K_F+M_{a_n+t}'|_F|\\
&\supset& |K_F+(b_n+t)\sigma^*(K_{F_0})|\\
&\supset& |(b_n+t+1)\sigma^*(K_{F_0})|.
\end{eqnarray*}

Take $t=p-1$. Noting that
$$|K_{X'}+M_{a_n+p-1}'+F|\subset |(a_n+p+m_0)K_{X'}|$$
and applying Lemma 2.7 of \cite{MPCPS} again, one has
$$a_{n+1}\pi^*(K_X)|_F\ge M_{a_n+p+m_0}|_F\ge M'_{a_n+p}|_F\ge
b_{n+1} \sigma^*(K_{F_0}).$$ Here we set $a_{n+1}:=a_n+p+m_0$ and
$b_{n+1}=b_n+p$. Set $\beta_n = \frac{b_{n}}{a_{n}}.$ Clearly
$\lim_{n\mapsto +\infty} \beta_n = \frac{p}{m_0+p}$.

The case $p=1$ can be proved similarly, but with a simpler
induction. We omit the details.
\end{proof}

The following Lemma is needed in our proof. Though similar one has
already been established in several papers of the first author, we
include it here for the convenience to readers.

\begin{lem}\label{b>0} Keep the same notation as in \ref{setup}.
Let $f:X'\longrightarrow B$ be the induced fibration of
$\varphi_{m_0}$. Denote by $F:=S$ a general fiber of $f$. If
$\dim(B)=1$ and $g(B)>0$, then $\pi^*(K_X)|_F\sim\sigma^*(K_{F_0})$
where $\sigma:F\longrightarrow F_0$ is the blow down onto the smooth
minimal model.
\end{lem}
\begin{proof} We shall use the idea of Lemma 14 in Kawamata's paper
\cite{KA}. By Shokurov's theorem in \cite{Sho}(see also \cite{HsM}),
each fiber of $\pi:X'\longrightarrow X$ is rationally chain
connected. Therefore, $f(\pi^{-1}(x))$ is a point for all $x\in X$.
Considering the image $G\subset (X \times B)$ of $X'$ via the
morphism $(\pi\times f)\circ \triangle_{X'}$ where $\triangle_{X'}$
is the diagonal map $X'\longrightarrow X'\times X'$, one knows that
$G$ is a projective variety. Let $g_1:G\longrightarrow X$ and
$g_2:G\longrightarrow B$ be two projections. Since $g_1$ is a
projective morphism and even a bijective map, $g_1$ must be both a
finite morphism of degree 1 and a birational morphism. Since $X$ is
normal, $g_1$ must be an isomorphism. So $f$ factors as $f_1 \circ
\pi$ where $f_1:=g_2\circ g_1^{-1} : X \rightarrow B$ is a well
defined morphisms. In particular, a general fiber $F_0$ of $f_1$
must be smooth minimal. So it is clear that
$\pi^*(K_X)|_F\sim\sigma^*(K_{F_0})$ where $\sigma$ is nothing but
$\pi|_F$.
\end{proof}

\section{\bf Proof of the main theorem}
We begin to prove Theorem \ref{main}. Let $X$ be a complex minimal
projective 3-fold of general type with $\mathbb{Q}$-factorial
terminal singularities and $\chi(\omega_X)\ge 0$.

\begin{setup}\label{reduction}{\bf Reduction to the case $p_g=1$.}
If $p_g(X)>1$, then $\varphi_8$ is birational by \cite{IJM}. If
$q(X)>0$, then $\varphi_m$ is birational for all $m\ge 7$ by
\cite{C-H}. Thus we assume, from now on, that $p_g(X)\le 1$ and
$q(X)=0$. The assumption
$$0\le \chi(\omega_X)=-\chi(\mathcal {O}_X)=-1+q(X)-h^2(\mathcal
{O}_X)+p_g(X)$$ implies $p_g(X)\ge h^2(\mathcal {O}_X)+1\ge 1$. Then
we clearly have $p_g(X)=1$, $h^2(\mathcal {O}_X)=0$ and
$\chi(\omega_X)=0$.
\end{setup}

According to \cite{Crelle}, we only have to study
non-Gorenstein minimal 3-folds $X$. In practice we may
assume the Cartier index $r(X)>1$.

\begin{setup}\label{plurigenus}{\bf Plurigenus.}
let $X$ be a minimal 3-fold pf general type with terminal
singularities. Recall Reid's plurigenus formula (at page
413 of \cite{YPG}):
$$P_m(X)=\frac{1}{12}m(m-1)(2m-1)K_X^3-(2m-1)\chi(\mathcal
{O}_X)+\sum_{Q} R_m(Q)  \eqno (4.1)$$ where $m>1$ is an integer, the
correction term
$$R_m(Q):=\frac{r^2-1}{12r}(m-\bar{m})+\sum_{j=0}^{\bar{m}-1}
\frac{\overline{bj}(r-\overline{bj})}{2r}$$ and the sum
$\sum_{Q}$ runs through all baskets of singularities $Q$ of
type $\frac{1}{r}(a,-a,1)$ with the positive integer $a$
coprime to $r$, $0<a<r$, $0<b<r$, $ab\equiv 1$ (\text{mod}
$r$) and $\bar{m}$ the smallest residue of m \text{mod}
$r$.  {\it Reid's result (Theorem 10.2 in \cite{YPG}) says
that the above baskets $\{Q\}$ of singularities are in fact
virtual (!) though $X$ may have  non-quotient terminal
singularities. Iano-Fletcher (\cite{Fletcher}) actually
shows that the set of baskets $\{Q\}$ in Reid's formula is
uniquely determined by $X$}.
\end{setup}

\begin{lem}\label{p3} For all basket $Q$, $R_5(Q)\ge R_4(Q)\ge
R_3(Q)\ge R_2(Q)$. In
particular,
$$P_5(X)>P_4(X)>P_3(X)>P_2(X)\ge 1$$ for all 3-fold $X$ with
$\chi(\mathcal {O}_X)=0$.
\end{lem}
\begin{proof} Suppose that $Q$ is of type
$\frac{1}{r}(a,-a,1)=\frac{1}{r}(1,-1, b)$ with $r>1$, $a$ coprime
to $r$, $ab\equiv 1$ (\text{mod} r) and $0<b<r$.

When $r=2$, one has
$$R_3(Q)=\frac{r^2-1}{6r}=\frac{1}{4}=R_2(Q).$$

When $r=3$, one has
$$R_3(Q)=\frac{r^2-1}{4r}=\frac{2}{3}>\frac{1}{3}=\frac{b(r-b)}{2r}=R_2(Q).$$

When $r>3$, one has $m=\bar{m}$ for $m=2,3$ and
$$R_3(Q)=\sum_{j=0}^2 \frac{\overline{bj}(r-\overline{bj})}{2r}\ge
\sum_{j=0}^1 \frac{\overline{bj}(r-\overline{bj})}{2r}=R_2(Q).$$

If $\chi(\mathcal {O}_X)=0$, then Reid's formula (4.1) gives
$$P_3(X)=\frac{5}{2}K_X^3+\sum_Q R_3(Q)>\frac{1}{2}K_X^3+\sum_Q R_2(Q)
=P_2(X)>0.$$
In particular, $P_3(X)\ge 2$.

Similarly one can verify the inequality $P_5(X)>P_4(X)>P_3(X)$. So
$P_5(X)\ge 4$.

\end{proof}

Now Theorem \ref{main} essentially follows from the
following theorem by letting $m_0=2, 3$.

\begin{thm}\label{m_0} Let $X$ be a minimal projective 3-fold of general type
with $\mathbb{Q}$-factorial terminal singularities, $p_g(X)=1$ and
$P_{m_0}(X)\ge 2$. Then $\varphi_{m}$ is birational onto its image
under one of the following conditions:

(1) $m\ge 4m_0+3$ and $2\le m_0\le 3$;

(2) $m\ge 4m_0+2$ and $m_0\ge 4$;

(3) $m\ge 14$, $\chi(\mathcal {O}_X)=0$ and $m_0=3$.
\end{thm}
\begin{proof}
Set $d_{m_0}:=\dim\varphi_{m_0}(X)$. We discuss according to the
value of $d_{m_0}$. We shall mainly apply Theorem \ref{Key}. Because
$p_g(X)>0$, Theorem \ref{Key}(i) is always satisfied whenever $m\ge
m_0+2$.

{\bf Case 1}. $d_{m_0}=3$. Recall that $S$ is a generic irreducible
element of $|M_{m_0}|$. We have $p=1$. On $S$ we take $G:=S|_S$.
Clearly $|G|$ is base  point free and is not composed of a pencil.
So a  generic irreducible element  $C$ of $|G|$ is a smooth curve.
Also under the assumption $m\ge m_0+2$ one has
\begin{eqnarray*}
&&K_S+\roundup{(m-1)\pi^*(K_X)-S-E_{m_0}'}|_S\\
&\ge &(K_{X'}|_S+\roundup{(m-1)\pi^*(K_X)-S-E_{m_0}'}|_S)+G\ge G.
\end{eqnarray*}

Thus Theorem \ref{Key}(ii) is also satisfied. Because
$m_0\pi^*(K_X)|_S\ge C$, we may take $\beta=\frac{1}{m_0}$ and so
Theorem \ref{Key}(iii) is satisfied.

Note that $C^2\ge 2$ because $|G|$ is not composed of a pencil. So
$m_0\pi^*(K_X)|_S\cdot C\ge C^2\ge 2$, which implies $\xi\ge
\frac{2}{m_0}$. Now we take $m\ge 3m_0+2$ and run Theorem \ref{Key}.
One has $\alpha=(m-1-\frac{m_0}{p}-\frac{1}{\beta})\xi\ge
2+\frac{2}{m_0}>2$.  This means that $\varphi_m$ is  birational for
all $m\ge 3m_0+2$. This is not the best. In fact, Theorem \ref{Key}
already gives $\xi\ge \frac{2g(C)+1}{3m_0+2}$. Note  that
$2g(C)-2=(K_S+C)\cdot C=(K_{X'}|_S+2C)\cdot C> 4$. One  has $\xi\ge
\frac{9}{3m_0+2}$. Now take $m>\frac{8}{3}m_0+\frac{13}{9}$. One has
$\alpha=(m-1-2m_0)\xi>2$. Theorem \ref{Key} says that $\varphi_m$ is
birational whenever $m>\frac{8}{3}m_0+\frac{13}{9}$. (One may go on
optimizing the estimate. We stop here since we have already proved
the theorem.)

{\bf Case 2}. $d_{m_0}=2$. We have $p=1$. On $S$ we take $G:=S|_S$.
Clearly $|G|$ is base  point free and is composed of a pencil. So a
generic irreducible element  $C$ of $|G|$ is a smooth curve. One has
$G\equiv tC$ for $t\ge P_{m_0}-2\ge 1$. Also under the assumption
$m\ge m_0+2$ one has
\begin{eqnarray*}
&&K_S+\roundup{(m-1)\pi^*(K_X)-S-E_{m_0}'}|_S\\
&\ge &(K_{X'}|_S+\roundup{(m-1)\pi^*(K_X)-S-E_{m_0}'}|_S)+G\ge C.
\end{eqnarray*}
So $|K_S+\roundup{(m-1)\pi^*(K_X)-S-E_{m_0}'}|_S|$ can separate
different generic irreducible elements of $|G|$ provided that $|G|$
is composed of a rational pencil. When $|G|$ is composed of an
irrational pencil, we need the assumption $m\ge 2m_0+2$. In fact, we
have $S|_S\equiv tC$ with $t\ge 2$ and $m_0\pi^*(K_X)|_S\equiv
tC+E_{m_0}'|_S$. Take two generic irreducible elements $C'$, $C''$
of $|G|$. Because
$$(m-m_0-1)\pi^*(K_X)|_S-C'-C''-\frac{2}{t}E_{m_0}'\equiv
(m-m_0-\frac{2}{t}m_0-1)\pi^*(K_X)|_S$$ is nef and  big, the
Kawamata-Viehweg vanishing theorem gives a surjective map
\begin{eqnarray*}
& &H^0(S,
K_S+\roundup{(m-m_0-1)\pi^*(K_X)|_S-\frac{2}{t}E_{m_0}'|_S})\\
&\longrightarrow & H^0(C', K_{C'}+D')\oplus H^0(C'', K_{C''}+D'')
\end{eqnarray*}
where $D'$, $D''$ are divisors of positive degree. Besides the last
groups are nonzero. Noting that
\begin{eqnarray*}
&&K_S+\roundup{(m-1)\pi^*(K_X)-S-E_{m_0}'}|_S\\
&\ge & K_S+\roundup{(m-m_0-1)\pi^*(K_X)|_S}\\
&\ge &K_S+\roundup{(m-m_0-1)\pi^*(K_X)|_S-\frac{2}{t}E_{m_0}'|_S},
\end{eqnarray*}
$|K_S+\roundup{(m-1)\pi^*(K_X)-S-E_{m_0}'}|_S|$ separates $C'$ and
$C''$. So we see that Theorem \ref{Key}(ii) is satisfied whenever
$m\ge 2m_0+2$.

Because $m_0\pi^*(K_X)|_S\ge C$, we may take $\beta=\frac{1}{m_0}$
and so Theorem \ref{Key}(iii) is satisfied.

If we take a big $m$ such that $\alpha$ is big enough, then Theorem
\ref{Key} gives
$$m\xi\ge 2g(C)-2+(m-1-\frac{m_0}{p}-\frac{1}{\beta}),$$
which says $\xi\ge \frac{2}{2m_0+1}$. This is only an initial
estimate. Take $m> 4m_0+2$. Then $\alpha=(m-2m_0-1)\xi>2$. Theorem
\ref{Key} says that $\varphi_m$ is  birational and  that $\xi\ge
\frac{5}{4m_0+3}$. Take $m=3m_0+2$. Then $\alpha=(m-2m_0-1)\xi>1$.
Theorem \ref{Key} gives $\xi\ge \frac{4}{3m_0+2}$. Finally take
$m>\frac{7}{2}m_0+2$. Then $\alpha=(m-2m_0-1)\xi>2$. Theorem
\ref{Key} says that $\varphi_m$ is  birational whenever
$m>\frac{7}{2}m_0+2$.

{\bf Case 3}. $d_{m_0}=1$. We have an induced fibration
$f:X'\longrightarrow B$. Denote by $F:=S$ a general fiber of $f$.
Note that $F$ is a surface of general type. Denote by
$\sigma:F\longrightarrow F_0$ the blowing down onto the minimal
model.

{\bf Subcase 3.1}. $g(B)>0$. By Lemma \ref{b>0}, we have
$\pi^*(K_X)|_F\sim \sigma^*(K_{F_0})$. We need to study the
condition (ii) in Theorem \ref{Key}. For all $m\ge m_0+5$ one has
\begin{eqnarray*}
&& K_F+\roundup{(m-1)\pi^*(K_X)-F-\frac{1}{p}E_{m_0}'}|_F\\
&\ge &K_F+\roundup{(m-m_0-1)\pi^*(K_X)|_F}\\
&\ge & K_F+(m-m_0-1)\pi^*(K_X)|_F\\
&\ge &(m-m_0)\pi^*(K_X)|_F\ge 5\sigma^*(K_{F_0}).
\end{eqnarray*}
So $\Phi_{|K_F+\roundup{(m-1)\pi^*(K_X)-F-\frac{1}{p}E_{m_0}'}|_F|}$
already gives a birational map. The proof of Theorem \ref{Key} tells
that we may omit the conditions $(iii)\sim (v)$ in Theorem
\ref{Key}. Thus $\varphi_m$ is birational for all $m\ge m_0+5$.

{\bf Subcase 3.2}. $g(B)=0$. For a general $F$ the natural map
$$H^0(X', K_{X'}-F)\longrightarrow H^0(X', K_{X'})$$
is a strict inclusion simply because $F$ is movable. This  means
that $p_g(F)>0$ since $p_g(X')=p_g(X)>0$. For our purpose we
classify $F$ into the following three types:
\begin{quote}
\item(I) $(K_{F_0}^2, p_g(F)=(1,2)$;

\item(II) $(K_{F_0}^2, p_g(F)=(2,3)$;

\item(III) $F$ does not belong to types (I) and (II).
\end{quote}
\medskip

We first study the type (III) case. By the results of Bombieri
\cite{Bom}, Reider \cite{Reider}, Catanese-Ciliberto \cite{C-C} and
P. Francia \cite{Francia}, $|2\sigma^*(K_{F_0})|$ is always base
point free whenever $p_g(F)>0$. We set $G:=2\sigma^*(K_{F_0})$ to
run Theorem \ref{Key}. The inclusion $\mathcal {O}(1)\hookrightarrow
f_*\omega_{X'}^{m_0}$ implies the inclusion
$$f_*\omega_{X'/\mathbb{P}^1}^2\hookrightarrow
f_*\omega_{X'}^{4m_0+2}.$$ Viehweg (\cite{VV}) first showed that
$f_*\omega_{X'/\mathbb{P}^1}^2$ is semi-positive. So it is also
generated by global sections. Thus it is  clear that
$$|M_{4m_0+2}||_F\supset |2\sigma^*(K_{F_0})|$$
where $M_{4m_0+2}$ is the movable part of $|(4m_0+2)K_{X'}|$.
\begin{quote}
(\#) So Theorem \ref{Key}(ii) is satisfied for all $m\ge 4m_0+2$ and
for all $F$ with $p_g(F)>0$, $G\le 2\sigma^*(K_{F_0})$ and $|G|$ not
an irrational pencil.
\end{quote}
In type (III) case, $G$ is an even divisor under our setting and is
not composed of a pencil of curves. On the other hand the
birationality of $\Phi_{|3K_F|}$ implies that a general $C\in |G|$
is non-hyperelliptic. By Lemma \ref{beta} we can take a
$\beta\mapsto \frac{1}{2m_0+2}$ such that $\pi^*(K_X)|_F-\beta C$ is
numerically equivalent to an effective $\mathbb{Q}$-divisor. By our
definition, one has $p\ge P_{m_0}-1\ge 1$. For all $m\ge 3m_0+4$,
$\alpha=(m-1-\frac{m_0}{p}-\frac{1}{\beta})\xi>0$. In general,
whenever $m\ge 4m_0+2$, Theorem \ref{Key} asserts that $\varphi_m$
is birational.

Next we study the type (II) case. We have  $(K_{F_0}^2,
p_g(F))=(2,3)$. By \cite{BPV} (page 227) we know that
$|\sigma^*(K_{F_0})|$ is base point free and $\Phi_{|K_F|}$ is
finite of degree 2. We set $G:=\sigma^*(K_{F_0})$. Then $C$, as a
generic irreducible element of $G$, is smooth and of  genus 3. We
have already showed (statement (\#) in type (III) case) that Theorem
\ref{Key}(ii) is satisfied for all $m\ge 4m_0+2$ because $G\le
2\sigma^*(K_{F_0})$. By Lemma \ref{beta} we can take a $\beta\mapsto
\frac{1}{m_0+1}$ such that $\pi^*(K_X)|_F-\beta C$ is numerically
equivalent to an effective $\mathbb{Q}$-divisor. So
$\xi=\pi^*(K_X)|_F\cdot C\ge \beta C^2\mapsto \frac{2}{m_0+1}$.
Taking the limit one has $\xi\ge \frac{2}{m_0+1}$. Similarly we have
$p\ge 1$. Take $m=3m_0+2$. Then
$\alpha=(m-1-\frac{m_0}{p}-\frac{1}{\beta})\xi\ge
\frac{2m_0}{m_0+1}>1$. Theorem \ref{Key} gives $\xi\ge
\frac{6}{3m_0+2}.$ In order to get the birationality we need the
assumption $m\ge 4m_0+2$, under which one has
$$\alpha=(m-1-\frac{m_0}{p}-\frac{1}{\beta})\xi\ge
\frac{12m_0}{3m_0+2}>2.$$ Thus $\varphi_m$ is  birational whenever
$m\ge 4m_0+2$.

Finally we study the type (I) case. We have $(K_{F_0}^2,
p_g(F))=(1,2)$. By \cite{BPV} we know that the movable part of
$|K_F|$ has one simple base point. Take $|G|$ to be the movable part
of $|K_F|$. Then a generic irreducible element $C$ of $|G|$ is a
smooth curve of genus 2. Similarly we need the assumption $m\ge
4m_0+2$ to secure the condition Theorem \ref{Key}(ii) (by virtue of
the statement (\#) in the type (III) case because $|G|$ is a
rational pencil). By Lemma \ref{beta} we can take a $\beta\mapsto
\frac{1}{m_0+1}$ such that $\pi^*(K_X)|_F-\beta C$ is numerically
equivalent to an effective $\mathbb{Q}$-divisor. Clearly we have
$p\ge 1$. If we take a big $m$ such that $\alpha$ is big enough,
Theorem \ref{Key} gives
$$m\xi\ge 2g(C)-2+(m-2m_0-2)\xi$$
which implies $\xi\ge \frac{1}{m_0+1}$. Take $m=4m_0+5$. Then
$\alpha\ge \frac{2m_0+3}{m_0+1}>2$. Theorem \ref{Key} gives $\xi\ge
\frac{5}{4m_0+5}$ and so does the birationality of $\varphi_m$. Take
$m=4m_0+4$. Then $\alpha\ge \frac{10m_0+10}{4m_0+5}>2$. One has the
birationality of $\varphi_m$ and $\xi\ge \frac{5}{4m_0+4}$. Take
$m=4m_0+3$. Then $\alpha\ge \frac{10m_0+5}{4m_0+4}>2$ whenever
$m_0\ge 2$. So $\varphi_{4m_0+3}$ is birational and $\xi\ge
\frac{5}{4m_0+3}$. We have already proved Theorem \ref{m_0}(1).

Assume $m_0\ge 4$ and take $m=4m_0+2$. Then $\alpha\ge
\frac{10m_0}{4m_0+3}>2$. Theorem \ref{Key} gives the birationality
of $\varphi_{4m_0+2}.$ We have proved Theorem \ref{m_0}(2).

The only left case to verify is $m_0=3$ and $\chi(\mathcal {O}_X)=0$
for which we have already the birationality of $\varphi_{4m_0+3}.$
We have
$$\xi\ge \frac{1}{3}\eqno (4.2)$$
as shown above.

\begin{quote}
(**)Furthermore the assumption $m_0=3$, $\chi(\mathcal {O}_X)=0$ and
$\xi>\frac{1}{3}$ gives $\alpha\ge 2m_0\xi>2$, which means that
$\varphi_{4m_0+2}=\varphi_{14}$ is birational by Theorem \ref{Key}.
\end{quote}

\begin{claim}\label{claim} When $m_0=3$ and $\chi(\mathcal {O}_X)=0$,
$\varphi_{4m_0+2}=\varphi_{14}$ is birational.
\end{claim}

To prove the claim, we need to study the 5-canonical system
$|5K_{X'}|$. By Lemma \ref{p3}, one has $P_5(X)\ge 4$. Denote
$d_5:=\dim\varphi_5(X)$. Keep the same notation as in \ref{setup}
and in Theorem \ref{Key}. Recall that $M_5$ is the movable part of
$|5K_{X'}|$. The induced fibration from $\varphi_5$ is denoted by
$f_5:X'\longrightarrow W_5$. (In fact we can take further
modifications over $X'$ (still denote by $X'$ the final
modification) such that the final $X'$ dominates both $\varphi_3$
and $\varphi_5$.)

If $d_5=3$, then $\dim \varphi_5(F)=2$ and $\dim\varphi_5(C)=1$ for
a general curve $C\in |G|$ in a  general fiber $F$. This means that
the linear system $|M_5||_C\ (\subset |M_5|_C|)$ gives a finite map
from $C$ onto a curve and so does $|M_5|_C|$. The Riemann-Roch and
Clifford theorem on $C$ says that $5\pi^*(K_X)\cdot C\ge M_5\cdot
C\ge 2$, i.e. $\xi\ge \frac{2}{5}>\frac{1}{3}$. Statement (**)
 tells that $\varphi_{14}$ is birational.

If $d_5=1$, then $|5K_{X'}|$ induces the same fibration
$f:X'\longrightarrow W_5=B$ simply because $5K_{X'}\ge 3K_{X'}$. The
typical property here is that $5\pi^*(K_X)\ge 3F$ for a general
fiber $F$. By Lemma \ref{beta} (just take $m_0'=5$ and $p'=3$) we
can find a $\beta\mapsto \frac{3}{8}$ such that $\pi^*(K_X)|_F-\beta
C$ is numerically equivalent to an effective $\mathbb{Q}$-divisor.
Going on the argument just before the claim, we have $m_0=3$, $p=1$,
$\beta=\frac{3}{8}$ and $\xi\ge \frac{1}{3}$ as in (4.2). Take
$m=14$. Then $\alpha=(m-1-m_0-\frac{1}{\beta})\xi>2$. So Theorem
\ref{Key} gives the birationality of $\varphi_{14}$.

Finally if $d_5=2$, we need to study the image surface $W_5'$ of
$X'$ through the morphism $\Phi_{|M_5|}$. In fact we have  the
decomposition $$\Phi_{|M_5|}: X' \overset{f_5} \longrightarrow W_5
\overset{s_5}\longrightarrow W_5'\subset \mathbb{P}^{P_5(X)-1}.$$
There is a very ample divisor $H_5$ on $W_5'$ such that
$M_5=\Phi_{|M_5|}^*(H_5)$. Furthermore one has $M_5|_{S_5}\equiv
a_5\tilde{C}$ for a general member $S_5\in |M_5|$ and an integer
$a_5\ge \deg{s_5}\deg(W_5')\ge \deg(W_5')\ge P_5(X)-2$, $\tilde{C}$
is a general fiber of $f_5$. If $a_5\ge 3$, we may utilize our
argument in {\bf Case 2} replacing $m_0$ by $m_0'=5$. Over there we
have shown $\xi\ge \frac{4}{3m_0'+2}=\frac{4}{17}$. But we can take
a better $\beta$, namely $\beta=\frac{3}{5}$. We have $p=1$. Take
$m=12$. Then $\alpha=(12-1-5-\frac{5}{3})\xi>1$. Theorem \ref{Key}
gives $\xi\ge \frac{1}{3}$. Take $m=14$. Then
$\alpha=(14-1-5-\frac{5}{3})\xi\ge \frac{19}{9}>2$. Theorem
\ref{Key} says that $\varphi_{14}$ is birational.

So we are left to study what happens if $a_5=2$. This means that
$\deg(W_5')=2$, $P_5(X)=4$ and  $\deg(s_5)=1$. Recall that a degree
2 irreducible surface in $\mathbb{P}^3$ must be one of the
following type of surfaces (see for instance Reid's lecture notes
\cite{Park}, Ex. 19 at page 30):

{\bf (a)} $W_5'$ is the cone $\overline{\mathbb{F}}_2$ by blowing
down the unique section with the self-intersection $-2$ on the
Hirzebruch surface $\mathbb{F}_2$ (a ruled surface);

{\bf (b)} $W_5'=\mathbb{P}^1\times \mathbb{P}^1$.
\medskip

We study these two cases separately in the following propositions.
\end{proof}

\begin{prop}\label{a} For case ${\bf (a)}$, $\varphi_{14}$ is
birational.
\end{prop}
\begin{proof}  We know that $M_5=g_5^*(H_5')$ for a very ample
divisor $H_5'$ on $W_5'$ with ${H_5'}^2=2$ and
$g_5:X'\longrightarrow W_5'$ is the birational morphism. Because
$W_5'$ is already normal (which is a cone), by taking further
modification to $X'$, we can assume that $g_5$ factors through the
minimal resolution $\mathbb{F}_2$ of $W_5'$. So we have the map
$g_5:X'\overset{h_5}\longrightarrow
\mathbb{F}_2\overset{\nu}\longrightarrow W_5'$ where $h_5$ is a
fibration and $\nu$ is the minimal resolution. Set
$\bar{H_5}=\nu^*(H_5')$. Then $\bar{H}^2= (H_5')^2 = 2$. Noting that
$\bar{H_5}$ is nef and big on $\mathbb{F}_2$, we can write
$$\bar{H}\sim \mu G_0+nT$$
where $G_0$ is the unique section with $G_0^2=-2$, $\mu$ and $n$ are
integers and $T$ is the general fiber of the ruling on
$\mathbb{F}_e$. The property of $\bar{H_5}$ being nef and big
implies $\mu>0$ and $n\ge 2\mu\ge 2$. Now let
$\theta_2:\mathbb{F}_2\longrightarrow \mathbb{P}^1$ be the ruling,
whose fibers are all smooth rational curves. Set $f_0:=\theta_2\circ
h_5: X'\longrightarrow \mathbb{P}^1$, which is a fibration with
connected fibers. Denote by $\tilde{F}$ a general fiber of $f_0$ and
by $\tilde{\sigma}:\tilde{F}\longrightarrow \tilde{F_0}$ the
contraction onto the minimal model. Clearly $p_g(\tilde{F})>0$. We
have
$$M_5\sim g_5^*(H_5')=h_5^*(\bar{H})\ge 2\tilde{F}.$$

Set $N_5=2\tilde{F}$. Replace $|M_5|$ by the sub-linear system
$|N_5|$ in Theorem \ref{Key}. We can also study the birationality of
$\varphi_m$ as remarked in \ref{remark}.

Note that the fibration $f_0:X'\longrightarrow \mathbb{P}^1$ is the
induced fibration by $\Phi_{|2\tilde{F}|}$. We can repeat a similar
argument to that in {\bf Subcase 3.2} above. We can verify those
conditions in Theorem \ref{Key}.

Because $14K_{X'}\ge 2\tilde{F}$, Theorem \ref{Key}(i) is satisfied.
On the other hand, we have $\mathcal {O}(2)\hookrightarrow
{f_0}_*\omega_{X'}^5$ which implies the inclusion
$${f_0}_*\omega_{X'/\mathbb{P}^1}^2\hookrightarrow
{f_0}_*\omega_{X'}^{12}.$$ Similarly Viehweg's semi-positivity and
the base point freeness of $|2\sigma^*(K_{\tilde{F_0}})|$ say that
$|12K_{X'}||_{\tilde{F}}$ can separate different generic irreducible
elements of $|2\sigma^*(K_{\tilde{F_0}})|$. So can
$|14K_{X'}||_{\tilde{F}}$. So Theorem \ref{Key}(ii) is satisfied.

Now we take $m_0'=5$ to run Theorem \ref{Key} with $|M_5|$ replaced
by $|N_5|$. One has $p=2$.

If $\tilde{F}$ is of type (III), we take
$\tilde{G}:=2\tilde{\sigma}^*(K_{\tilde{F_0}})$. Take $m=14$. The
proof of Lemma \ref{beta} says that one can take a $\beta\mapsto
\frac{1}{7}$ such that $\pi^*(K_X)|_{\tilde{F}}-2\beta
\tilde{\sigma}^*(K_{\tilde{F_0}})$ is numerically equivalent to an
effective $\mathbb{Q}$-divisor. So
$\alpha=(m-1-\frac{m_0'}{p}-\frac{1}{\beta})\tilde{\xi}>0$. Noting
that a generic irreducible element of $|\tilde{G}|$ is
non-hyperelliptic and even, $\varphi_{14}$ is birational by Theorem
\ref{Key}.

If $\tilde{F}$ is of type (II), we set
$\tilde{G}:=\tilde{\sigma}^*(K_{\tilde{F_0}})$. Then $\tilde{C}$, as
a generic irreducible element of $|\tilde{G}|$, is smooth and of
genus 3. We have already showed above that Theorem \ref{Key}(ii) is
satisfied for $m=14$. By Lemma \ref{beta} we can take a
$\beta\mapsto \frac{2}{7}$ such that $\pi^*(K_X)|_{\tilde{F}}-\beta
\tilde{C}$ is numerically equivalent to an effective
$\mathbb{Q}$-divisor. So $\tilde{\xi}=\pi^*(K_X)|_{\tilde{F}}\cdot
\tilde{C}\ge \beta \tilde{C}^2\mapsto \frac{4}{7}$. Taking the limit
one has $\tilde{\xi}\ge \frac{4}{7}$. We have $p=2$. Take $m=14$.
Then $\alpha=(m-1-\frac{m_0'}{p}-\frac{1}{\beta})\tilde{\xi}\ge 4$.
Theorem \ref{Key} says that $\varphi_{14}$ is birational.

If $\tilde{F}$ is of type (I), we take $\tilde{G}$ to be the movable
part of $|\tilde{\sigma}^*(K_{\tilde{F_0}})|$. Then $\tilde{C}$, as
a generic irreducible element of $|\tilde{G}|$, is smooth and of
genus 2. We have already showed above that Theorem \ref{Key}(ii) is
satisfied for $m=14$. By Lemma \ref{beta} we can take
$\beta=\frac{2}{7}$ such that $\pi^*(K_X)|_{\tilde{F}}-\beta
\tilde{\sigma}^*(K_{\tilde{F_0}})$ is numerically equivalent to an
effective $\mathbb{Q}$-divisor. So
$\tilde{\xi}=\pi^*(K_X)|_{\tilde{F}}\cdot C\ge \beta= \frac{2}{7}$.
Take $m=12$. Then
$\alpha=(m-1-\frac{m_0'}{p}-\frac{1}{\beta})\tilde{\xi}>1$. Theorem
\ref{Key} gives $\tilde{\xi}\ge \frac{1}{3}$. Take $m=14$. Then
$\alpha=(m-1-\frac{m_0'}{p}-\frac{1}{\beta})\tilde{\xi}\ge
\frac{7}{3}>2$. Theorem \ref{Key} says that $\varphi_{14}$ is
birational.
\end{proof}

\begin{prop}\label{b} For case ${\bf (b)}$, $\varphi_{14}$ is
birational.
\end{prop}
\begin{proof} Recall that we have  two fibrations
$f:X'\longrightarrow B$ and $f_5:X'\longrightarrow
W_5'=\mathbb{P}^1\times \mathbb{P}^1$. For a general canonical curve
$C$ in the a general fiber $F$ of $f$, we can study
$\dim\varphi_5(C)$. If $\dim\varphi_5(C)=1$, then $\Phi_{|M_5|}$
maps $C$ onto a curve. Clearly $5\pi^*(K_X)\cdot C\ge M_5\cdot C\ge
2$. So $\xi\ge \frac{2}{5}>\frac{1}{3}$ and $\varphi_{14}$ is
birational according to the statement (**). {}From now on we assume
that $\dim\varphi_5(C)=0$ for a general $C$. We may take further
blowing ups to $X'$ so that a general $C$ is simply a generic fiber
of $f_5$.

Because the only very ample divisor $H_5'$ on $W_5'$ with
${H_5'}^2=2$ is $L_1+L_2=q_1^*(\text{point})+q_2^*(\text{point})$
where $q_1, q_2$ are projection maps from $\mathbb{P}^1\times
\mathbb{P}^1$ to $\mathbb{P}^1$. Set $\tilde{f_i}:=q_i\circ f_5:
X'\longrightarrow \mathbb{P}^1$, $i=1,2$. Then $\tilde{f_1}$ and
$\tilde{f_2}$ are two fibrations onto $\mathbb{P}^1$. Let $F_1$ and
$F_2$ are respectively general fibers of $\tilde{f_1}$ and
$\tilde{f_2}$. Then $F_1\cap F_2$ is simply a general fiber $C$ of
$f_5$. We will prove alternatively that $\varphi_{14}$ is
birational.

Consider the sub-liner system
$|K_{X'}+\roundup{8\pi^*(K_X)}+F_1+F_2|\subset |14K_{X'}|$ which
clearly separates different fibers of $\tilde{f_1}$. Take a general
$F_1$ as a fiber of $\tilde{f_1}$. Because $8\pi^*(K_X)+F_2$ is nef
and big, the Kawamata-Viehweg vanishing theorem gives a surjective
map
$$H^0(K_{X'}+\roundup{8\pi^*(K_X)}+F_2+F_1)\longrightarrow
H^0(F_1, K_{F_1}+\roundup{8\pi^*(K_X)}|_{F_1}+C).$$ We hope to prove
that $|K_{F_1}+\roundup{8\pi^*(K_X)}|_{F_1}+C|$ gives a birational
map for a general $F_1$. Note that $|C|$ is a rational pencil on the
surface $|F_1|$. Clearly $|K_{F_1}+\roundup{8\pi^*(K_X)}|_{F_1}+C|$
separate different $C$ simply because $p_g(F_1)>0$ and
$\roundup{8\pi^*(K_X)}|_{F_1}\ge 0$. Take a general curve $C$ in the
family $|F_2|_{F_1}|$. The vanishing theorem again gives the
surjective map
$$H^0(F_1, K_{F_1}+\roundup{8\pi^*(K_X)|_{F_1}}+C)\longrightarrow
H^0(C, K_C+D)$$ where $D:=(\roundup{8\pi^*(K_X)|_{F_1}})|_C$ is a
divisor of $$\deg(D)\ge (8\pi^*(K_X)|_{F_1})\cdot C=8\pi^*(K_X)\cdot
C=8\xi\ge \frac{8}{3}>2$$ recalling the inequality (4.2) and that
$C$ is also a canonical curve in the surface $F$. So
$\Phi_{|K_C+D|}$ is an embedding and the birationality principle
says that $\varphi_{14}$ is birational. We have completely showed
Theorem \ref{m_0}.
\end{proof}

\begin{exmp}\label{examp2} Example \ref{examp1} shows that Theorem
\ref{m_0} is optimal for $m_0=3$. It is also optimal for $m_0=2$. In
fact, Fletcher has an example (page 151, No. 18 in \cite{C-R}): the
canonical hyper-surface $X_{22}\subset \mathbb{P}(1,2,3,4,11).$ It
is clear $m_0=2$ and that $\varphi_{11}$ is birational, but
$\varphi_{10}$ is not birational.

Recently we were informed of a new construction of a 3-fold $Y$ by
E. Stagnaro (\cite{S}). $Y$ also has $p_g=1$, $P_2=2$ and
$\varphi_{10}$ is not birational.
\end{exmp}
\bigskip

By a similar method we are able to prove Corollary \ref{corollary}.
To avoid unnecessary redundancy we omit some details.

\begin{setup}\label{sketch}{\bf A sketch proof of Corollary
\ref{corollary}}
\begin{proof}
First one has $P_2(X)\ge 4$ by a similar application of Reid's
formula (see the proof of Lemma \ref{p3}). So one may take $m_0=2$.
Set $d_2:=\dim \varphi_2(X)$.

If $d_2\ge 2$, it is clear that Theorem \ref{Key}(i), (ii) are
satisfied for $m\ge 8$. Now the situation $d_2=3$ follows directly
from the argument in {\bf Case 1} of Theorem \ref{m_0}. In fact, one
get the birationality of $\varphi_{7}$. In the situation $d_2=2$,
one may also follow the argument in {\bf Case 2} of Theorem
\ref{m_0}. The only difference is that we can take $p=2$ here. So it
can be verified that $\varphi_8$ is birational.

If $d_2=1$, one has an induced fibration $f:X'\longrightarrow B$.
When $g(B)>0$, the argument in {\bf Subcase 3.1} of Theorem
\ref{m_0} shows that $\varphi_{m}$ is birational for all $m\ge 7$.
We assume $g(B)=0$ from now on. Let $F$ be a general fiber of $f$.
If $p_g(F)=0$, then $q(X)=q(F)=0$ and the assumption
$\chi(\omega_X)>0$ implies $p_g(X)\ge 2$. So $\varphi_{m}$ is
birational for all $m\ge 8$ by the main theorem in \cite{IJM}. We
are reduced to study the situation $p_g(F)>0$. Now we can follow the
argument of {\bf Subcase 3.2} in Theorem \ref{m_0}. The typical
property here is the inclusion
$$\mathcal {O}(3)\hookrightarrow f_*\omega_{X'}^2.$$
We can still choose $G\le 2\sigma^*(K_{F_0})$ on $F$. We can verify
that Theorem \ref{Key} (i), (ii) are satisfied for all $m\ge 6$.
Also we can choose much better $p$ ($p\ge 2$) and $\beta$,
respectively, than those in {\bf Subcase 3.2} in Theorem \ref{m_0}.
We can verify with less difficulties, case by case (for type (I),
(II), (III)) by applying Theorem \ref{Key}, that $\varphi_{m}$ is
birational for all $m\ge 8$.

In a word, Corollary \ref{corollary} is true.
\end{proof}
\end{setup}

\begin{setup}{\bf Final remarks.} Clearly, for the case $p_g=0$,
Theorem \ref{Key} will generate some new results which
improves a corollary of J. Koll\'ar (\cite{Kol}) and the
first author on the birationality of $\varphi_m$, where $m$
is certain function in terms of $m_0$. We do not include
our result here because we are not sure whether that would
be optimal.
\end{setup}


\end{document}